\magnification \magstep1
\def \qed {\vrule height6pt  width6pt depth0pt}
\centerline {\bf Operators which factor through Banach lattices 
not containing $c_0$}

\medskip

\centerline {by}

\medskip

\centerline {N. Ghoussoub and W.~B. 
Johnson\footnote*{Supported in part by NSF DMS--8703815}}

\bigskip

In this supplement to [GJ1], [GJ3], we give an intrinsic 
characterization of (bounded, linear) operators on Banach 
lattices which factor through Banach lattices not containing a 
copy of $c_0$ which complements the characterization of [GJ1], 
[GJ3] that an operator admits such a factorization if and only if 
it can be written as the product of two operators neither of 
which preserves a copy of $c_0$.  The intrinsic characterization 
is that the restriction of the second adjoint of the operator to 
the ideal generated by the lattice in its bidual does not 
preserve a copy of $c_0$.  This property of an operator was 
introduced by C. Niculescu [N2] under the name ``strong type 
B".  That this condition is stronger than the property that the 
operator itself does not preserve a copy of $c_0$ is exhibited 
by the example constructed in [FGJ] and studied further in 
[GJ2].

When we wrote [GJ1], [GJ3] we did not know that Niculescu 
[N1], [N2], [N3] was independently considering some of the 
same problems; we now want to acknowledge his priority for 
discovering some results that overlap with ours.  Also, in 
remark 1 below, we indicate how our results, when translated 
into Niculescu's language, answer some questions he raised.  

Before stating the theorem of this note, we introduce some 
notation: Given a Banach lattice $E$, $I(E)$ (respectively, 
$Bd(E)$) denotes the closed order ideal (respectively, band) in 
$E^{**}$ generated by $E$.  The symbols $E$, $F$, $G$ will be 
reserved for Banach lattices while $X$, $Y$, $Z$ will be used for 
general Banach spaces.  Terminology is as in [GJ3], and the 
reader will probably need a copy of [GJ3] on hand to follow the 
proof of the theorem.

{\proclaim Theorem.  Let $E$ be a Banach lattice and let $T$ be 
an operator from $E$ into a Banach space $X$.  The following 
are equivalent:

(1) ${T^{**}}_{|I(E)}$ does not preserve a copy of $c_0$.

(2) ${T^{**}}_{|I(E)}$ does not preserve a positive disjoint copy 
of $c_0$.

(3) ${T^{**}}[Bd(E)] \subset X$.

(4) There exists a Banach space $Y$ and bounded linear 
operators $u\colon E \to Y$ and $v\colon Y \to X$ such that 
neither $u$ nor $v$ preserves a copy of $c_0$ and $T=vu$.

(5) There exists a weakly sequentially complete Banach lattice 
$F$, an interval preserving lattice homomorphism $R\colon I(E) 
\to F$, and an operator $S\colon F \to X$ so that 
${T^{**}}_{|I(E)}=SR$.}

\medskip
\noindent {\bf Proof.} That (1) and (2) are equivalent follows 
from Theorem I.3 of [GJ3] (that the other conditions in Theorem 
I.3 are equivalent to (2) was proved earlier by Niculescu---see 
Theorem 3.4 in [N2]).  Niculescu [N2, Proposition 4.2] stated the 
equivalence of (2) and (3); we follow his lead in leaving the 
proof as an exercise for the reader, but note that (2) is the 
statement we tie to the other conditions. The equivalence of (4) 
with (5) was proved in Corollary I.8 in [GJ3] and (5) 
$\Rightarrow$ (2) is immediate, so we only need to check that 
(2) implies (5).  In fact, the proof of this implication is 
essentially contained in the proof of Theorem I.7 of [GJ3], but 
since only weaker implications are explicitly stated there we 
present here a detailed proof.  

Let $W$ be the absolute convex solid hull of $T^*Ball(X^*)$ in 
$E^*$ and let $G$ be the interpolation space (cf. [GJ3]) 
corresponding to $W$ in $E^*$.  Let $J \colon G \to E^*$ be the 
canonical injection and $I \colon X^* \to G$ such that $T^* = JI$.  
By Lemma I.5 of [GJ3] we have for each $\mu$ in $E^{**}$ and 
$n=1,2,\dots,$
$$\Vert J^*\mu \Vert \le 2^n sup \,  \{ \, \Vert T^{**} \nu 
\Vert \colon \, \nu \in
E^{**}, |\nu | \le |\mu | \, \} + 2^{-n} \Vert \mu \Vert. \eqno 
(\dag )$$

Letting $F$ be the norm closure of $J^*[I(E)]$ in $G^*$ and 
setting  $R={J^*}_{|I(E)}$,  $S = {I^*}_{|F}$, we have that 
${T^{**}}_{|I(E)}$ factors through $R$ and hence in order to 
complete the proof we only need to check that $F$ does not 
contain a copy of $c_0$.

Since $I(E)$ is an ideal in $E^{**}$ we have from $(\dag )$ that 
for each $f \in I(E)$ and $n=1,2,\dots ,$
$$\Vert R f \Vert \le 2^n sup \,  \{ \, \Vert T^{**} g \Vert 
\colon \, g \in I(E),
|g | \le |f | \,\} + 2^{-n} \Vert f \Vert. \eqno (\ddag )$$

\noindent  Now $J^*$ is a lattice homomorphism since $J$ is 
interval preserving ([AB, p. 90]).  Since ${T^{**}}_{|I(E)}$ does 
not preserve a copy of $c_0$, neither does $R$ by Lemma I.1.a 
of [GJ3], hence (\ddag ) implies that $F$ has order continuous 
norm by Lemma I.1.b of [GJ3].

We first prove that if $(e_m)_{m=1}^{\infty}$ is an increasing 
sequence in $I(E)_+$ such that $(J^* e_m) = (R e_m)$ is norm 
bounded in $F$ then it is norm convergent in $F$.  Indeed, let $ 
\mu $ be the $w^*$--limit of $(J^* e_m)_{m=1}^{\infty}$ in 
$G^*$.  Since $J^*[E^{**}]$ is norm dense in $G^*$, there exists for 
each $\epsilon > 0$ an element $\nu_{\epsilon} \in E^{**}$ such 
that $\Vert J^* \nu_{\epsilon} - \mu \Vert \le \epsilon$. 

Since $I(E)$ is an ideal in $E^{**}$, the sequence 
$(\nu_{\epsilon} \wedge e_m)_{m=1}^{\infty}$ is in $I(E)$; 
moreover, it is increasing and norm bounded by $\Vert 
\nu_{\epsilon}  \Vert$, hence $J^*(\nu_{\epsilon} \wedge 
e_m)$ norm converges to some $z_{\epsilon}$ in $F$ by 
Theorem I.3.d in [GJ3].  Using again the fact that $J^*$ is a 
lattice homomorphism, we have that $J^*(\nu_{\epsilon} 
\wedge e_m) = J^*(\nu_{\epsilon}) \wedge J^*( e_m) $ hence 
$J^*(\nu_{\epsilon} \wedge e_m)$ $w^*$--converges in $G^*$ to 
$J^*(\nu_{\epsilon}) \wedge \mu$.  Therefore 
$J^*(\nu_{\epsilon}) \wedge \mu = z_{\epsilon}$ is in $F$.  
Thus we get 
$$\lim_{\epsilon \to 0} \Vert z_{\epsilon} - \mu \Vert = 
\lim_{\epsilon \to 0} \Vert (J^* \nu_{\epsilon} )  \wedge \mu - 
\mu \Vert \le \lim_{\epsilon \to 0} \Vert J^* \nu_{\epsilon} - 
\mu \Vert = 0.$$

\noindent It follows that $\mu$ is in $F$.  Since the sequence 
$(J^*e_m)_{m=1}^{\infty}$ is in the order interval $[0, \mu]$ 
and $F$ is order continuous, it follows that 
$(J^*e_m)_{m=1}^{\infty}$ norm converges to $\mu $. 

Suppose now that $(f_m)_{m=1}^{\infty}$ is any increasing 
sequence in $F$.  To show that it is norm convergent, it is 
enough to show that for any $\epsilon > 0$ there exists 
$(e_m)_{m=1}^{\infty}$ increasing in $I(E)_+$ such that $\Vert 
f_m - J^*e_m \Vert \le \epsilon $.  For that, let $z_m \in I(E)_+$ 
be such that  $\Vert J^*z_m - f_m \Vert \le \epsilon 2^{-m}$.  
Set $e_m = \vee_{i=1}^m z_i$.  Then $(e_m)_{m=1}^{\infty}$ is 
increasing and 
$$|J^* e_m - f_m | = \left | J^*\left ( \vee_{i=1}^m z_i \right ) - 
f_m  \right | = \left | \vee_{i=1}^m J^* z_i - \vee_{i=1}^m f_i 
\right | \le \sum_{i=1}^m | J^* z_i - f_i|.$$
\noindent It follows that
$$\Vert J^* e_m - f_m \Vert \le \sum_{i=1}^m \Vert f_i - J^* z_i 
\Vert \le \epsilon \sum_{i=1}^m 2^{-i} \le \epsilon.\eqno $$

By the classical version of Theorem I.3 in [GJ3], $F$ does not 
contain a copy of $c_0$ and hence $F$ is a weakly sequentially 
complete Banach lattice ([LT, p. 34]). \hfill \qed

\medskip

In the above proof of (5) $\Rightarrow $ (2), if $X$ is Banach 
lattice and $T$ is a positive operator, then every positive 
operator dominated by $T$ also factors through $R$.  Thus we 
get:
\medskip

\proclaim Corollary.  Let $E$ be a Banach lattice, let $T$ be a 
positive operator from $E$ into a Banach lattice  $X$, and 
assume that $ {T^{**}}_{|I(E)}$ does not preserve a positive 
disjoint copy of $c_0$.  Then there exists a weakly sequentially 
complete Banach lattice $F$ and an interval preserving lattice 
homomorphism $R\colon I(E) \to F$ such that for any operator 
$ 0 \le T_1 \le T$ there is a positive operator $S_1\colon F \to 
X$ so that ${T^{**}_1}_{|I(E)}=S_1 R$.
\bigskip

\noindent {\bf Remark 1.} The examples in [FGJ] of positive 
quotient mappings onto $c_0$ which do not preserve a copy of 
$c_0$ were shown in [GJ2] to fail property (4) of the theorem, 
hence are, in Niculescu's terminology [N2], operators of type B 
which are not of strong type B.  This answers problem 4.1 in 
[N2].  (However, if $E$ is order continuous then $I(E)=E$ and so 
operators from $E$ of type B are of strong type B.)  These 
operators also give a solution to problem 4.4 in [N2]; namely, 
there exists a lattice homomorphism $T$ on some Banach 
lattice which does not preserve a copy of $c_0$ or of $l_1$, yet 
$T^2$ is not weakly compact, while $T^3$ must be weakly 
compact: see remark 2 in [GJ3, p. 164].  
\medskip

\noindent {\bf Remark 2.} One of the operators mentioned in 
remark 1 (``$T_1$" in the terminology of [GJ2]) maps weakly 
Cauchy sequences into weakly convergent sequences; this gives 
a counterexample to one part of Proposition 4.2 in [N2].
\medskip

\noindent {\bf Remark 3.} Since a weakly sequentially 
complete Banach lattice has the so-called {\it analytic Radon-
Nikodym---ARN---property} (cf. [BD]), an operator which 
satisfies condition (3) of the theorem is an ARN-operator.  This 
implication was recently proved in a direct way by Bukhvalov 
[B].  Although aware of [GJ3], Bukhvalov mistakenly thought 
that his result could not follow from the factorizations in [GJ3] 
(he applied the error from [N2] mentioned in remark 2 to 
conclude that the operators mentioned in remark 1 satisfy 
condition (3) of the theorem). 

\bigskip

\centerline {\bf References}

\bigskip

\item{[AB]} C.~D. Aliprantis and O. Burkinshaw, {\it Positive 
operators,} {\bf Academic Press} (1984).
\smallskip

\item {[BD]} A.~V. Bukhvalov and A.~A. Danilevich, {\it 
Boundary properties of analytic and harmonic functions with 
values in Banach space,} {\bf Mat. Zametki 31} (1982), 203--
214.  English translation: {\bf Math. Notes 31} (1982), 104--
110.
\smallskip

\item {[B]} A.~V. Bukhvalov, {\it On the analytic Radon-
Nikodym property,} [unpublished manuscript].
\smallskip

\item {[FGJ]} T. Figiel, N. Ghoussoub, and W.~B. Johnson, {\it On 
the structure of non-weakly compact operators on Banach 
lattices,} {\bf Math. Ann. 257} (1981), 317--334.
\smallskip

\item {[GJ1]} N. Ghoussoub and W.~B. Johnson, {\it On subspaces 
of Banach lattices not containing $C(\Delta )$,} [unpublished].
\smallskip

\item {[GJ2]} N. Ghoussoub and W.~B. Johnson, {\it 
Counterexamples to several problems on the factorization of 
bounded linear operators,} {\bf Proc. AMS 92} (1984), 233--
238.
\smallskip

\item {[GJ3]} N. Ghoussoub and W.~B. Johnson, {\it Factoring 
operators through Banach lattices not containing $C(0,1)$,} {\bf 
Math. Zeit. 194} (1987), 153--171.
\smallskip

\item {[LT]} J. Lindenstrauss and L. Tzafriri, {\it Classical 
Banach spaces II: Function spaces,} {\bf Ergebnisse Math. 
Grenzgebiete 97} Springer-Verlag (1979).
\smallskip

\item {[N1]} C.~Niculescu, {\it Weak compactness in Banach 
lattices,} {\bf J. Operator Theory 6} (1981), 217--231.
\smallskip

\item {[N2]} C.~Niculescu, {\it Order $\sigma$-continuous 
operators on Banach lattices,} {\bf Lecture Notes in Math. 991} 
(1983), 188--201.
\smallskip

\item {[N3]} C.~Niculescu, {\it Operators of Type A and local 
absolute continuity,} {\bf J. Operator Theory 13} (1985), 49--
61.
\bigskip

\noindent Department of Mathematics, University of British 
Columbia, Vancouver, B.~C., Canada
\smallskip
\noindent Department of Mathematics, Texas A\&M University, 
College Station, TX 77843, USA

\vfill

\end